\documentclass[11pt]{amsart}

\newcommand{\eps}{\epsilon}

\newfont{\fnt}{cmr10 scaled 550}
\renewcommand{\eps}{\varepsilon}

\newtheorem{conj}{Conjecture}
\newtheorem{lemma}{Lemma}
\newtheorem{cor}{Corollary}
\newtheorem{prop}{Proposition}

\theoremstyle{remark}
\newtheorem{remark}{Remark}

\font\strange=msbm10

\renewcommand{\epsilon}{\varepsilon}

\renewcommand{\Sigma}{\varSigma}
\newcommand{\R}{{{\mathchoice  {\hbox{$\textstyle{\text{\strange R}}$}}
{\hbox{$\textstyle{\text{\strange R}}$}}
{\hbox{$\scriptstyle  N\kern-0.3em  R$}}  
{\hbox{$\scriptscriptstyle  R\kern-0.2em  R$}}}}}

\newcommand{\Z}{{{\mathchoice  {\hbox{$\textstyle{\text{\strange Z}}$}}
{\hbox{$\textstyle{\text{\strange Z}}$}}
{\hbox{$\scriptstyle  Z\kern-0.3em  Z$}}
{\hbox{$\scriptscriptstyle  Z\kern-0.2em  Z$}}}}}

\newcommand{\N}{{{\mathchoice  {\hbox{$\textstyle{\text{\strange N}}$}}
{\hbox{$\textstyle{\text{\strange N}}$}}
{\hbox{$\scriptstyle  N\kern-0.3em  N$}}
{\hbox{$\scriptscriptstyle  N\kern-0.2em  N$}}}}}

\renewcommand{\phi}{\varphi}
\usepackage{amsmath,amsthm,amscd}

\begin{document}
\title{Proof of the normal scalar curvature conjecture}
\date{August 26, 2007}

 \author{Zhiqin Lu}

\address{Department of
Mathematics, University of California,
Irvine, Irvine, CA 92697}


\email[Zhiqin Lu]{zlu@math.uci.edu}

\thanks{The
 author is partially supported by  NSF Career award DMS-0347033 and the
Alfred P. Sloan Research Fellowship.}

\maketitle 

\pagestyle{myheadings}

\newcommand{\ka}{K\"ahler }
\newcommand{\ii}{\sqrt{-1}}
\section{Introduction}
Let $M^n$ be an  $n$-dimensional manifold isometrically immersed into the space form $N^{n+m}(c)$ of constant sectional curvature $c$. Define the normalized scalar curvature $\rho$ (resp. $\rho^\perp$) for the tangent bundle (resp.  the normal bundle)  as follows:
\begin{align}
\begin{split}
&\qquad \quad\quad \rho=\frac{2}{n(n-1)}\sum_{1=i<j}^n R(e_i, e_j,e_j,e_i),\\
&(resp.)\quad \rho^{\perp}=\frac{2}{n(n-1)}\left(\sum_{1=i<j}^n\sum_{1=r<s}^m\langle
R^\perp(e_i,e_j)\xi_r,\xi_s\rangle^2\right)^{\frac 12},
\end{split}
\end{align}
where $\{e_1,\cdots, e_n\}$ (resp. $\{\xi_1,\cdots,\xi_m\}$) is an orthonormal basis of the tangent (resp. normal) bundle, and $R$ (resp. $R^\perp$) is the curvature tensor for the tangent (resp. normal) bundle.

In the study of submanifold theory, De Smet, Dillen, Verstraelen, and Vrancken~\cite{ddvv} made the following {\sl normal scalar curvature conjecture}\footnote{Also known as the DDVV conjecture.}: 

\begin{conj}\label{cc1} Let $h$ be the second fundamental form, and 
let $H=\frac 1n \,{\rm trace}\, h$ be the mean curvature tensor.  Then 
\[
\rho+\rho^\perp\leq |H|^2+c.
\]
\end{conj}

 Let $x\in M$ be a fixed point and let $(h_{ij}^r)$ ($i,j=1,\cdots,n$ and $r=1,\cdots,m$) be the entries  of (the traceless part of ) the second fundamental form under the  orthonormal bases of both the tangent bundle and the normal bundle. Then by
 ~\cite{bogd}, or ~\cite{dfv}, Conjecture~\ref{cc1} can be formulated as an inequality with respect to the coefficients $(h_{ij}^r)$ as follows:
  \begin{align}\label{1a}
\begin{split}
&\sum_{r=1}^m\sum_{1=i<j}^n(h_{ii}^r-h_{jj}^r)^2+2n\sum_{r=1}^m\sum_{1=i<j}^n(h_{ij}^r)^2\\
&\geq 2n\left(\sum_{1=r<s}^m\sum_{1=i<j}^n\left(\sum_{k=1}^n(h_{ik}^rh_{jk}^s-h_{ik}^sh_{jk}^r)\right)^2\right)^{\frac 12}.
\end{split}
\end{align}

Suppose that $A_1,A_2,\cdots,A_m$ are $n\times n$ symmetric real matrices. Let
\[
||A||^2=\sum_{i,j=1}^n a_{ij}^2,
\]
where $(a_{ij})$ are the entries of $A$, 
and let
\[
[A,B]=AB-BA
\]
be the commutator. Then the inequality~\eqref{1a}, in terms of matrix notations, can be formulated as
\begin{conj}\label{conj2}
 For $n, m\geq 2$, we have
\begin{equation}\label{conj}
(\sum_{r=1}^m||A_r||^2)^2\geq 2(\sum_{r<s}||[A_r,A_s]||^2).
\end{equation}
\end{conj}

Fixing $n,m$, we call the above inequality Conjecture $P(n,m)$. Note that Conjecture~\ref{cc1} is equivalent to Conjecture~\ref{conj2}, which is purely  linear algebraic.

 A weaker version of  Conjecture~\ref{cc1},
$
 \rho\leq |H|^2+c,
 $
 was proved in ~\cite{ch1}. An alternate proof is in~\cite{bogd2}.

The following special cases of Conjecture~\ref{conj2} were known. $P(2,m)$ and $P(n,2)$ were proved in~\cite{ddvv}; $P(3,m)$ was proved in~\cite{cl}; and $P(n,3)$ was proved in~\cite{lu21}, the previous version of this paper. In~\cite{dfv}, a weaker version of $P(n,m)$ was proved by using an algebraic inequality in~\cite{lianmin} (see also~\cite{xusenlin}) . In the same paper, $P(n,m)$ was proved under the addition assumption that the submanifold is either {\it Lagrangian $H$-umbilical}, or {\it ultra-minimal} in $\mathbb C^4$.

 In this  paper, we prove the conjecture for any $n,m\geq 1$.

\section{Invariance}

Let $A_1,\cdots,A_m$ be  $n\times n$ symmetric matrices.
Let $G=O(n)\times O(m)$. Then $G$ acts on matrices $(A_1,\cdots,A_m)$ in the following natural way: let $(p,q)\in G$, where $p,q$ are $n\times n$ and $m\times m$ orthogonal matrices, respectively. Let $q=\{q_{ij}\}$. Then
\[
(p,I)\cdot (A_1,\cdots,A_m)=(pA_1p^{-1},\cdots, pA_mp^{-1}),
\]
and 
\[
(I,q)\cdot (A_1,\cdots,A_m)=(\sum_{j=1}^m q_{1j}A_j,\cdots,\sum_{j=1}^mq_{mj}A_{j}).
\]

It is easy to verify the following 

\begin{prop}\label{prop1}
Conjecture $P(n,m)$ is $G$ invariant. That is, in order to prove inequality ~\eqref{conj}
for $(A_1,\cdots,A_m)$, we just need to prove the inequality for any $\gamma\cdot (A_1,\cdots,A_m)$ where $\gamma\in G$. Moreover, the expressions of both sides of ~\eqref{conj} are $G$ invariant.
\end{prop}

\qed

\begin{cor}
We can prove Conjecture~\ref{conj2} under the following additional assumptions on the matrices:
\begin{enumerate}
\item $A_1$ is diagonal;
\item $\langle A_\alpha,A_\beta\rangle=0$ if $\alpha\neq \beta$;
\item $||A_1||\geq\cdots\geq ||A_m||$.
\end{enumerate}
\end{cor}

Note that under the above assumptions, $A_k=0$ if $k>\frac 12n(n+1)$.

\section{Proof of $P(n,m)$.}

In this section, we prove Conjecture~\ref{conj2}. We first establish some lemmas which are themselves interesting.

\begin{lemma} Suppose $\eta_1\cdots,\eta_n$ are real numbers and 
\[
\eta_1+\cdots+\eta_n=0,\quad
\eta_1^2+\cdots+\eta_n^2=1.
\]
Let $r_{ij}\geq 0$ be nonnegative numbers for $i<j$. Then we have
\begin{equation}\label{use}
\sum_{i<j}(\eta_i-\eta_j)^2 r_{ij}\leq\sum_{i<j}r_{ij}+{\rm Max} (r_{ij}).
\end{equation}
\end{lemma}

{\bf Proof.}
We assume that
$
\eta_1\geq\cdots\geq\eta_n.
$
If $\eta_1-\eta_n\leq 1$ or $n=2$, then ~\eqref{use} is trivial. So we assume $n>2$,  and
\[
\eta_1-\eta_n>1.
\]
We observe that $\eta_i-\eta_j\leq 1$ for $2\leq i<j\leq n-1$. Otherwise, we could have
\[
1\geq\eta_1^2+\eta_n^2+\eta_i^2+\eta_j^2>\frac 12
((\eta_1-\eta_n)^2+(\eta_i-\eta_j)^2)>1,
\]
which is a contradiction. 

Using the same reason,
if  $\eta_1-\eta_{n-1}>1$,  then we have $\eta_2-\eta_{n}\leq 1$; and if $\eta_2-\eta_{n}>1$, then we have $\eta_1-\eta_{n-1}\leq 1$.
Replacing $\eta_1,\cdots,\eta_n$ by $-\eta_n\cdots,-\eta_1$ if necessary, we can always assume that $\eta_2-\eta_{n}\leq 1$.
Thus $\eta_i-\eta_j\leq 1$ if $2\leq i<j$, and ~\eqref{use} is implied by the following inequality
\begin{equation}\label{use2}
\sum_{1<j}(\eta_1-\eta_j)^2 r_{1j}\leq\sum_{1<j} r_{1j}+\underset{1<j}{\rm Max}\, (r_{1j}).
\end{equation}
Let $s_j=r_{1j}$ for $j=2,\cdots,n$. 
The the above inequality becomes
\begin{equation}\label{use3}
\sum_{1<j}(\eta_1-\eta_j)^2 s_{j}\leq\sum_{1<j} s_{j}+\underset{1<j}{\rm Max}\, (s_{j}).
\end{equation}

In order to prove the above inequality, we define the matrix $P$ as follows
\[
P=\begin{pmatrix}
\underset{{1<j}}{\sum} s_j& -s_{2}&\cdots &-s_{n}\\
-s_{2}& s_{2}\\
\vdots&&\ddots
\\
-s_{n}&&&s_{n}
\end{pmatrix}.
\]

We 
claim that the maximum eigenvalue of $A$ is no more than $r=\sum_j s_j+{\rm Max}\, (s_j)$. To see this, we compute the determinant of the matrix
\[
\begin{pmatrix}
y-\underset{{1<j}}{\sum} s_j& s_{2}&\cdots &s_{n}\\
s_{2}& y-s_{2}\\
\vdots&&\ddots
\\
s_{n}&&&y-s_{n}
\end{pmatrix}.
\]
Using the Cramer's rule,  the answer is
\[
(y-s_2)\cdots (y-s_n)\left( y-\sum_{1<j}s_j-\sum_{1<j}\frac{s_j^2}{y-s_j}\right).
\]
For any $y>r$, we have $y-s_j>\sum_{s=2}^n s_j$. Thus the above expression is greater than
\[
(y-s_2)\cdots (y-s_n)(y-\sum_{1<j}s_j-(\sum_{1<j}s_j)^{-1}\sum_{1<j} s_j^2)>0.
\]

Let $\eta=(\eta_1,\cdots,\eta_n)^T$, we then have
\[
\sum_{1<j}(\eta_1-\eta_j)^2 s_{j}=
\eta^TP\eta\leq r=\sum_{1<j}^m s_j+\underset{1<j}{\rm Max}\, (s_j).
\]

\qed

\begin{lemma}\label{lem2}
Let $A$ be an $n\times n$  diagonal  matrix of norm $1$. Let $A_2,\cdots, A_m$ be symmetric matrices such that
\begin{enumerate}
\item $\langle A_\alpha,A_\beta\rangle=0$ if  $\alpha\neq \beta$;
\item $||A_2||\geq\cdots\geq||A_m||$.
\end{enumerate}
Then we have
\begin{equation}\label{key-1}
\sum_{\alpha=2}^m||[A,A_\alpha]||^2\leq \sum_{\alpha=2}^m ||A_\alpha||^2+||A_2||^2.
\end{equation}
\end{lemma}

{\bf Proof.} Replacing each $A_\alpha$ with its off-diagonal component won't change the left hand side of the above inequality, but will decrease the right hand side of the above.
 Thus without loss of generality, we  assume that each $A_\alpha$ has zero diagonal component.  Furthermore, we assume that each $A_\alpha$ is not zero.

Let $A_\alpha=((a_\alpha)_{ij})$, where  $(a_\alpha)_{ij}$ are  the entries for $\alpha=2,\cdots,m$. Let
\[
\delta=\underset{i\neq j}{\rm Max}\,\sum_{\alpha=2}^m(a_\alpha)^2_{ij}.
\]
Let 
\[
A=\begin{pmatrix}
\eta_1\\
&\ddots\\
&&\eta_n
\end{pmatrix}.
\]
Then by the previous lemma, we have
\begin{equation}\label{ab-1}
\sum_{\alpha=2}^m||[A,A_\alpha]||^2\leq \sum_{\alpha=2}^m||A_\alpha||^2+2\delta.
\end{equation}
Thus it remains to prove that
\begin{equation}\label{fun}
2\delta\leq ||A_2||^2.
\end{equation}
To see this, we identify each $A_\alpha$ with the (column)   vector $\vec A_\alpha$ in $\mathbb R^{\frac 12 n(n+1)}$. Let $\mu_\alpha$ be the norm of the vector. Then we have
\begin{equation}\label{cd-1}
\mu_\alpha^2=\frac 12 ||A_\alpha||^2
\end{equation}
for $\alpha=2,\cdots,m$. Extending the set of vectors $\{\vec A_\alpha/\mu_\alpha\}_{2\leq\alpha\leq m}$ into an orthonormal basis of $\mathbb R^{\frac 12n(n+1)}$
\[
\vec A_2/\mu_2, \cdots, \vec A_m/\mu_m,\vec A_{m+1},\cdots, \vec A_{\frac 12n(n+1)+1},
\]
we get an orthogonal matrix. Apparently, each row vector of the matrix is a unit vector. Thus we have
\[
\sum_{\alpha=2}^m (\mu_\alpha)^{-2}(a_\alpha)^2_{ij}\leq 1.
\]
Since $\mu_2\geq\cdots\geq\mu_m$, we get
\[
\sum_{\alpha=2}^m (a_\alpha)^2_{ij}\leq\mu_2^2\leq \frac 12||A_2||^2.
\]
This proves ~\eqref{fun}.

\qed

\begin{remark} Let $A$ be a diagonal matrix of unit norm and let $B$ be a symmetric matrix. Let $||B||_\infty={\rm Max}\, (|b_{ij}|)$, where $(b_{ij})$ are the entries of $B$. By~\eqref{ab-1}, we get
\[
||[A,B]||^2\leq ||B||^2+2||B||_\infty^2.
\]
Although not used directly in this paper, this is the crucial step of estimate that makes the whole proof work. Note that in~\cite{ddvv} (or by $P(n,2)$), we have  a much weaker version of the above inequality
\[
||[A,B]||^2\leq 2||B||^2.
\]
\end{remark}

{\bf Proof of Conjecture~\ref{conj2}.} 
Let $a>0$ be the largest positive real number such that
\[
(\sum_{\alpha=1}^m||A_\alpha||^2)^2\geq 2a(\sum_{\alpha<\beta}||[A_\alpha,A_\beta]||^2).
\]
Since $a$ is maximum, 
by the invariance, we can find matrices $A_1,\cdots,A_m$ such that
\begin{equation}\label{conj-1}
(\sum_{\alpha=1}^m||A_\alpha||^2)^2= 2a(\sum_{\alpha<\beta}||[A_\alpha,A_\beta]||^2)
\end{equation}
with the following additional properties:
\begin{enumerate}
\item $A_1$ is diagonal;
\item $\langle A_\alpha,A_\beta\rangle=0$ if $\alpha\neq \beta$;
\item $||A_1||\geq ||A_2||\geq\cdots \geq ||A_m||$.
\end{enumerate}

We let $t^2=||A_1||^2$ and let $A'=A_1/|t|$. Then ~\eqref{conj-1} becomes a quadratic expression in terms of $t^2$:
\begin{align*}&
t^4-2t^2(a\sum_{1<\alpha}||[A',A_\alpha]||^2-\sum_{1<\alpha}||A_\alpha||^2)
+(\sum_{\alpha=2}^m||A_\alpha||^2)^2\\&\qquad- 2a(\sum_{1<\alpha<\beta}||[A_\alpha,A_\beta]||^2)=0.
\end{align*}
Since the left hand side of the above is non-negative for all $t^2$, we have
\[
a\sum_{1<\alpha}||[A',A_\alpha]||^2-\sum_{1<\alpha}||A_\alpha||^2>0,
\]
and 
\[
||A_1||^2=a\sum_{1<\alpha}||[A',A_\alpha]||^2-\sum_{1<\alpha}||A_\alpha||^2.
\]
By  Lemma~\ref{lem2}, we have
\[
\sum_{1<\alpha}||[A',A_\alpha]||^2\leq \sum_{\alpha=2}^m||A_\alpha||^2+||A_2||^2\leq \sum_{\alpha=1}^m||A_\alpha||^2,
\]
which proves that $a\geq 1$.

\qed

\section{Appendix: on the B\"ottcher-Wenzel Conjecture}
In the study of random matrices, B\"ottcher and Wenzel~\cite{bw} posed the  following conjecture:

\begin{conj}\label{bw-conj}
Let $X,Y$ be two $n\times n$ matrices. Then
\[
||[X,Y]||^2\leq 2||X||^2||Y||^2,
\]
where  the norm is defined as 
\[
||X||^2=\sum_{i,j=1}^n x_{ij}^2.
\]
\end{conj}
 B\"ottcher and Wenzel proved the following special cases of the conjecture: if $n=2$, or $X$ is of rank $1$, or $X$ is normal, then the conjecture is true. Furthermore, they proved the following weaker version of the conjecture:
 \[
 ||[X,Y]||^2\leq 3||X||^2||Y||^2.
 \]

\smallskip

\smallskip

In this Appendix, we prove Conjecture~\ref{bw-conj}.

 We fix $X$ and assume that $||X||=1$. Let $V=gl(n,\mathbb R)$. Define a linear map 
\[
T: V\to V,\quad Y\mapsto [X^T,[X,Y]].
\]
Then we have

\begin{lemma}\label{lem1}
$T$ is a semi-positive definite symmetric linear transformation of $V$.
\end{lemma}

{\bf Proof.} This is a straightforward computation
\[
\langle Y_1,[X^T,[X,Y_2]]\rangle=\langle[X,Y_1],[X,Y_2]\rangle
=\langle [X^T,[X,Y_1]],Y_2\rangle.
\]

Obviously $T$ is semi-positive.

\qed

The conjecture is equivalent to the statement that the maximum eigenvalue of $T$ is not more than  $2$.

We let $\alpha$ be the maximum eigenvalue of $T$. Then $\alpha>0$. Let $Y$ be an eigenvector of $T$ with respect to $\alpha$. Then we have
\[
T(Y)=\alpha Y.
\]
A straightforward computation gives
\[
T([X^T,Y^T])=\alpha [X^T,Y^T],
\]
where $X^T$ is the transpose of $X$.

We claim that $Y$ and $Y_1=[X^T,Y^T]$ are linearly  independent: first, $Y_1\neq 0$, and second $\langle Y,Y_1\rangle=0$. Thus, we have the  following conclusion
\begin{prop}
\label{prop2}
The multiplicity of the eigenvalue $\alpha$ is at least $2$.
\end{prop}
\qed

Let 
\[
X=Q_1\Lambda Q_2
\]
be the singular decomposition of $X$, where $Q_1,Q_2$ are orthogonal matrices and $\Lambda$ is a diagonal matrix. Let
\[
B=Q_2YQ_2^{-1},\quad C=Q_1^{-1}YQ_1.
\]
Then we have
\[
||[X,Y]||^2=||\Lambda B-C\Lambda||^2.
\]
Let
\[
\Lambda=\begin{pmatrix}
s_1\\&\ddots\\&&s_n
\end{pmatrix}.
\]
Without loss of generality, we assume that $s_1\geq\cdots \geq s_n$. Since $||X||=1$, we have
\[
s_1^2+\cdots+s_n^2=1.
\]

Assume that $s_1^2\leq 1/2$. Then we have
\begin{equation}\label{bw-4}
||\Lambda B-C\Lambda||^2=\sum_{i,j=1}^n(s_ib_{ij}-s_j c_{ij})^2
\leq \sum_{i,j=1}^n 2(b_{ij}^2+c_{ij}^2)s_1^2\leq 2.
\end{equation}

Thus in this case, the conjecture is trivially true.
Now assume that $s_1^2>1/2$. By Proposition~\ref{prop2}, we can find an eigenvector $Y$ of $T$ such that 1). $||Y||=1$, and 2). $b_{11}=0$.

The conjecture can be proved if we can prove that 
\[
||[X,Y]||^2\leq 2.
\]

We first have the following equality (because $b_{11}=0$)
\[
||\Lambda B-C\Lambda||^2=c_{11}^2 s_1^2
+\sum_{i=2}^n(s_ib_{i1}-s_1c_{i1})^2
+\sum_{j=2}^n (s_1b_{1j}-s_jc_{1j})^2+\Delta_1,
\]
where we define 
\[
\Delta=\sum_{i=2}^n b_{1i}^2+\sum_{i=1}^n c_{i1}^2,
\]
and 
\[
\Delta_1=\sum_{i,j=2}^n(s_ib_{ij}-s_j c_{ij})^2.
\]
Apparently we have
\[
\Delta_1\leq \sum_{i,j=2}^n (b_{ij}^2+c_{ij}^2),
\]
because $s_2^2\leq 1/2$.
Thus we just need to prove that 
\[
c_{11}^2s_1^2
+\sum_{i=2}^n(s_ib_{i1}-s_1c_{i1})^2
+\sum_{i=2}^n (s_1b_{1i}-s_jc_{1i})^2
\leq \Delta+\sum_{i=2}^n b_{i1}^2+\sum_{i=2}^n c_{1i}^2.
\]

We consider the matrix
\[
P=
\begin{pmatrix}
\Delta& -b_{12}c_{12}-b_{21}c_{21}&\cdots &-b_{1n}c_{1n}-b_{n1}c_{n1}\\
-b_{12}c_{12}-b_{21}c_{21}& b_{21}^2+c_{12}^2\\
\vdots&&\ddots\\
-b_{1n}c_{1n}-b_{n1}c_{n1}&&&b_{n1}^2+c_{1n}^2
\end{pmatrix}.
\]

The above inequality is equivalent to that
 the maximum eigenvalue of the above matrix is no more than $\Delta+\sum_{i=2}^n b_{i1}^2+\sum_{i=2}^n c_{1i}^2$. To see this, we let
\[
y=\Delta+\sum_{i=2}^n b_{i1}^2+\sum_{i=2}^n c_{1i}^2+\eps
\]
for $\eps>0$. We have
\[
\det(yI-P)=\prod_{i=2}^n (y-b_{i1}^2-c_{1i}^2)\left(
y-\Delta-\sum_{i=2}^n \frac{(b_{1i}c_{1i}+b_{i1}c_{i1})^2}{y-b_{i1}^2-c_{1i}^2}\right).
\]
Let
\[
\beta=\max (b_{i1}^2+c_{1i}^2).
\]
Then we have
\[
y-\Delta-\sum_{i=2}^n \frac{(b_{1i}c_{1i}+b_{i1}c_{i1})^2}{y-b_{i1}^2-c_{1i}^2}\geq \beta+\eps-\beta\sum_{i=2}^n\frac{b_{1i}^2+c_{i1}^2}{\sum_{i=2}^n( b_{1i}^2+c_{i1}^2)}>0.
\]
The conjecture is proved.

\qed

\bibliographystyle{abbrv}   
\bibliography{pub,unp,2007}   
\vspace{0.2in}
Add in the proof: Recently J. Ge and Z. Tang gave an independent proof of the normal scalar curvature conjecture.
 
\end{document}